\documentstyle[twocolumn]{article}
\makeatletter
\def\@normalsize{\@setsize\normalsize{12pt}\xpt\@xpt
\abovedisplayskip 10pt plus2pt minus5pt\belowdisplayskip \abovedisplayskip
\abovedisplayshortskip \z@ plus3pt\belowdisplayshortskip 6pt plus3pt
minus3pt\let\@listi\@listI}

\def\subsize{\@setsize\subsize{12pt}\xipt\@xipt}

\def\section{\@startsection {section}{1}{\z@}{24pt plus 2pt minus 2pt}
{12pt plus 2pt minus 2pt}{\large\bf}}

\def\subsection{\@startsection {subsection}{2}{\z@}{12pt plus 2pt minus 2pt}
{12pt plus 2pt minus 2pt}{\subsize\bf}}
\makeatother

\newcommand{\pf}{\noindent{\bf proof\ }}
\newcommand{\qed}{\penalty 1000 \hfill\penalty 1000$\Box$\par\medskip}
\newenvironment{proof}{\pf }{\qed }
\newcounter{myfig}
\newcommand{\startfig}{\begin{figure}\refstepcounter{Theorem}
                      \setcounter{myfig}{\value{Theorem}} }
\newcommand{\eJJ}{\end{figure}}
     
\newcounter{myeqn}

\newcommand{\fo}{\hbox{{\rm FO}}}
\newcommand{\lfp}{\mbox{{\rm LFP}}}
\newcommand{\cownt}{\hbox{{\rm COUNT}}}
\newcommand\G{{\cal G}}        
\newcommand{\var}{\hbox{{\rm var}}}
\newcommand{\free}{\hbox{{\rm free}}}
\newcommand\D{{\cal D}}        
\newcommand{\bigset}[2]{\bigl\{ #1 \,\bigm|\, #2 \bigr\} }
\newcommand{\abs}[1]{ \vert #1 \vert }
\renewcommand\H{{\cal H}}        
\newcommand{\esub}{\preceq}        
\newcommand\sel{{\cal L}}      
\newcommand{\th}{^{\rm th}}     
\newcommand\B{{\cal B}}        
\newcommand{\bit}{\hbox{{\rm BIT}}}
\newcommand{\ac}{\hbox{{\rm AC}}}
\newcommand{\apal}{{\it Annals of Pure and Applied Logic }}
\renewcommand{\angle}[1]{ \langle #1 \rangle }
\newcommand{\pspace}{\hbox{{\rm PSPACE}}}
\newcommand{\p}{\hbox{{\rm P}}}
\newcommand{\set}[1]{\{  #1 \} }
\newcommand{\cq}[2]{(\exists #1\,#2)}    
\newcommand\J{{\cal J}}        
\newcommand{\iterp}{\hbox{{\rm ITER}}}
\newcommand\C{{\cal C}}        

\newcommand\nc{\newcommand} \nc\rnc{\renewcommand} 
\def\mboxes(#1){\xmboxes#1,xxx,}\def\endpiece{xxx}
\def\xmboxes#1,{\def\tmp{#1}%
\ifx\tmp\endpiece \else\expandafter\def\csname #1\endcsname{\hbox{#1}} 
\expandafter\xmboxes \fi}

\mboxes(Dom,HE,MP,MS,OSet,Proj,Range,Set,Vocabulary)

\nc\nn{\newenvironment} \nc\nt{\newtheorem}
\nt{Theorem}{Theorem}[section]
\nt{Corollary}[Theorem]{Corollary}
\nt{DefinitionAux}[Theorem]{Definition}
\nn{Definition}{\begin{DefinitionAux}\em}{\qed\end{DefinitionAux}}
\nt{Fact}[Theorem]{Fact}
\nt{Lemma}[Theorem]{Lemma}
\nt{Question}[Theorem]{Question}
\nt{ProvisoAux}{Proviso} 
\nn{Proviso}{\begin{ProvisoAux}\em}{\qed\end{ProvisoAux}}
\nt{RemarkAux}{Remark} 
\nn{Remark}{\begin{RemarkAux}\em}{\ $\Box$\end{RemarkAux}}

\nc\al{\alpha} \nc\be{\beta} \nc\ga{\gamma} \nc\Ga{\Gamma}
\nc\de{\delta}\nc\De{\Delta} \nc\ep{\varepsilon}
\rnc\phi{\varphi} \nc\si{\sigma}\nc\Si{\Sigma}
\nc\om{\omega}\nc\Om{\Omega}

\nc\cE{{\cal E}} \nc\cF{{\cal F}} \nc\cL{{\cal L}}
\rnc\L{L_{\infty,\omega}^\omega}
\nc\Lk{L_{\infty,\omega}^k}

\rnc\a{\bar{a}} \rnc\b{\bar{b}} \nc\n{\bar{n}} 
\rnc\u{\bar{u}} \rnc\v{\bar{v}} 
\nc\x{\bar{x}} \nc\y{\bar{y}} \nc\X{\bar{X}}

\nc\ns{\normalsize}
\nc\bs{\bigskip} \nc\ms{\medskip} \rnc\ss{\smallskip}
\rnc\ni{\noindent} \nc\ds{\displaystyle}
\nc\openmultiset{\{\!\!\{} \nc\closemultiset{\}\!\!\}}

\begin{document}

\date{}

\title{\Large\bf McColm's Conjecture}

\author{\begin{tabular}[t]{c@{\extracolsep{4em}}c@{\extracolsep{4em}}c}
  Yuri Gurevich\thanks{Partially supported by NSF, ONR and BSF.
Electrical Engineering and Computer Science Department, University of
Michigan, Ann Arbor, MI 48109-2122, USA, gurevich@umich.edu}	& 
Neil Immerman\thanks{Supported by NSF grant CCR-9207797.  Computer Science
Department, University of Massachusetts, Amherst, MA 01003, USA,
immerman@cs.umass.edu}  &
Saharon Shelah\thanks{Publication 525.  Partially supported by USA--Israel
Binational Science Foundation.  Mathematics Department, Hebrew University,
Jerusalem 91904, Israel, shelah@cs.huji.ac.il, and Mathematics Department,
Rutgers University, New Brunswick, NJ 08903, USA, shelah@math.rutgers.edu}\\
 \\
  U of Michigan & U of Massachusetts & Hebrew U and Rutgers\\
\end{tabular}}

\maketitle

\thispagestyle{empty}

\subsection*{\centering Abstract}
{\em Gregory McColm conjectured that positive elementary inductions are bounded in a class $K$ of
finite structures if every $(\fo + \lfp)$ formula is equivalent to a first-order formula in $K$.
Here $(\fo + \lfp)$ is the extension of first-order logic with the least fixed point operator.  We
disprove the conjecture.  Our main results are two model-theoretic constructions, one deterministic
and the other randomized, each of which refutes McColm's conjecture.}

\section{Introduction}

Gregory McColm conjectured in \cite{McColm} that, for every class $K$ of
finite structures, the following three claims are equivalent:
\begin{description}
\item[M1 ]
Every positive elementary induction is bounded in $K$.
\item[M2 ]
Every $(\fo + \lfp)$ formula is equivalent to a first-order formula in $K$.
\item[M3 ]
Every $L_{\infty\om}^\om$-formula is equivalent to a first-order formula in
$K$.
\end{description}

The definitions of $\L$ and $(\fo +\lfp)$ are recalled in the next section.

Clearly, M1 implies M2.  McColm observed that M3 implies M1.  Phokion Kolaitis and Moshe Vardi
proved that M1 implies M3 \cite{KV}.  A nice exposition of all of that is found in \cite{dawar}
The question whether M2 implies M1 has been open though McColm made the following important
observation.

Let $\n$ be the set $\{0,1,..,n-1\}$ with the standard order.  It is easy to
see that no infinite class of structures $\n$ satisfies M1.  List all
$(\fo + \lfp)$ sentences in vocabulary $\{<\}$: $\phi_0,\phi_1,\ldots$.  Let
\( K_i = \{\n\mid \n\models\phi_i\} \)
and construct an infinite $K$ such that every intersection $K\cap K_i$ is
either finite or co-finite.  Each $\phi_i$ is equivalent to a first-order
sentence in $K$.  Thus M1 does not follow from the restriction of M2 to
formulas without free variables.

The main results of this paper are two model theoretic constructions, one deterministic and the
other randomized, each of which gives a counterexample to the implication M2$\to$M1.  Actually, each
construction implies the stronger result that M2 fails to imply M1 even when $(\fo +\lfp)$ is
replaced in M2 by an arbitrary countable subset of $\L$, see Corollary \ref{L co} and Theorem
\ref{TheoremOne}.  We present the deterministic construction
in full detail in Section \ref{deterministic se}.  The randomized construction is presented in
Section \ref{random se}; but, some of the proofs are omitted due to lack of space.

Both constructions depend on the fact that the language $\L$, and thus $(\fo + \lfp)$ is unable to count the number
of vertices in a large clique.  The deterministic construction extends naturally to Theorem \ref{count th}:
an extension of our counterexample to the stronger language $(\fo +\lfp + \cownt)$ in which counting
is present.

Recall that $(\fo + \iterp)$, is first-order logic plus an
unbounded iteration operator (equivalent to the ``while'', and ``partial fixed point'' operators).
It is known that the language $(\fo + \iterp)$ captures $\pspace$ on ordered structures \cite{bounds,Vardi}.
Abiteboul and Vianu \cite{AV} showed that $\p =\pspace$ if and only if,  $(\fo + \lfp) = (\fo + \iterp)$ on all
sets of finite structures.  

In light of this, another interesting consequence of the deterministic construction is Corollary
\ref{ppspace co} which says that if P is not equal to PSPACE, then there is a set of finite
structures on which $\fo = (\fo +\lfp)$, but on which $\fo \ne (\fo + \iterp)$.

\section{Background}

We briefly recall some background material.  More information on Descriptive Complexity and Finite
Model Theory can be
found for example in \cite{ams} and \cite{Gurevich}.

\begin{Proviso}
Structures are finite.  Vocabularies are finite and do not contain function
symbols of positive arity.  In particular, the vocabulary of any
$\L$-formula is finite.  Classes of structures are closed under isomorphism.
\end{Proviso}

If $M$ is a structure then $|M|$ is the universe of $M$.  If $X$ is a
nonempty subset of $M$ (that is, of $|M|$) then $M\mid X$ is the induced
substructure with universe $X$. 

An {\em $r$-ary global relation\/} $\rho$ on a class $K$ of structures of
the same vocabulary is a function that, given a structure $M\in K$, produces
an $r$-ary (local) relation $\rho^M$ on $|M|$.  By definition,
$M\models\rho(\a)$ if and only if $\a\in\rho^M$.  It is supposed that, for
every isomorphism $\eta$ from $M$ to a structure $N$ and every $r$-tuple
$x_1,\ldots,x_r$ of elements of $M$, 
\( M\models\rho(x_1,\ldots,x_r) \iff N\models\rho(\eta(x_1,\ldots,x_r) )\).

In this paper, an {\em infinitary formula} means an
$L_{\infty,\omega}^\omega$ formula of finite vocabulary.  Recall that
$L_{\infty,\omega}^\omega$ is the generalization of first-order logic that
allows arbitrary infinite conjunctions and disjunctions provided that the
total number of individual variables, bound or free, in the resulting
formula is finite \cite{Barwise}.  In other words, infinitary formulas are
built from atomic formulas by means of negation, existential quantification,
universal quantification and the following rule:
\begin{itemize}\item
If $\{\phi_i\mid i\in I\}$ is a collection of infinitary formulas that uses
only a finite vocabulary and a finite number of individual variables
then
\( \bigvee_i\phi_i \mbox{\ \ \ and\ \ \ } \bigwedge_i\phi_i \)
are infinitary formulas.
\end{itemize}

The semantics is obvious.  $A\models\bigvee_i\phi_i(\a)$ if and only if
$A\models\phi_i(\a)$ for some $i$, and $A\models\bigwedge_i\phi_i(\a)$ if
and only if $A\models\phi_i(\a)$ for all $i$.  Let $L_{\infty,\omega}^k$ be the subset of
$L_{\infty,\omega}^\omega$ in which at most the $k$ distinct variables $\set{x_1,x_2,\ldots,x_k}$
occur. 

We next recall the definition of $(\fo +\lfp)$.
Consider a first-order formula
\( \phi(P,v_1,\dots,v_r,v_{r+1},\ldots,v_s) \)
with free individual variables $v_1,\ldots,v_s$ where an $r$-ary predicate
$P$ has only positive occurrences; let $\tau=\Vocabulary(\phi)-\{P\}$.  Given
a $\tau$-structure $M$ and elements $a_{r+1},\dots,a_s$ of $M$, we have the
following $r$-ary relations on the universe $|M|$ of $M$:
\[
P_0 \;=\; \emptyset,\qquad P_{i+1} \;=\]
\[ \{(v_1,\ldots,v_r)\mid M\models\phi(P_i,v_1,\ldots,v_r,a_{r+1},\ldots,a_s)\}\]
Since $P$ is positive in $\phi$, $P_0\subseteq P_1\subseteq
P_2\subseteq\ldots$.  M1 asserts that, for every such $\phi$, there exists a
positive integer $j$ such that, for every $M\in K$ and any
$a_{r+1},\dots,a_s\in M$, $P_j=\bigcup_i P_i$.

The least fixed point operator LFP can be applied to the formula $\phi$.
The result is a new formula
\[ \lfp_{P;v_1,\ldots,v_r}\phi(v_1,\ldots,v_s) \]
of vocabulary $\tau$.  If $M$ is a $\tau$-structure, $a_1,\ldots,a_s$ are
elements of $M$ and relations $P_i$ are as above then
\[M\models \lfp_{P;v_1,\ldots,v_r}\phi(a_1,\ldots,a_s)
\,\Leftrightarrow\, (a_1,\ldots,a_r)\in\bigcup_i P_i. \]
$(\fo + \lfp)$ is the extension of first-order logic with this new
formula-constructor.  Applications of LFP can be nested and interleaved with
other formula-constructors.  It is obvious that $(\fo + \lfp)$ is a subset of $\L$.

Pebble games are a convenient tool to deal with infinitary formulas.  A $k$-pebble game
$\Ga_\tau^k(A,B)$ is played by Spoiler and Duplicator on structures $A,B$ of vocabulary $\tau$.  For
each $i\in\{1,\ldots,k\}$, there are two pebbles numbered $i$; there are $2k$ pebbles altogether.
Starting with Spoiler, the players alternate making moves.  A move consists of placing a free pebble
at an element of one of the two structures or removing one of the pebbles from some element.  If
Spoiler puts a pebble of number $i$ at an element $x$ of $A$ (resp., an element $y$ of $B$),
Duplicator must answer by placing the other pebble number $i$ at some element $y$ of $B$ (resp.,
some element $x$ of $A$).  If Spoiler removes a pebble number $i$, Duplicator must remove the other
pebble number $i$.  Initially, all pebbles are free.  At each even-numbered state $S$, the pebbles
define a partial map $\eta_S$ from $A$ to $B$.  $\Dom(\eta_S)$ consists of the elements of $A$
covered by pebbles.  If $x\in A$ is covered by a pebble $i$ then $\eta_S(x)$ is the element of $B$
covered by the other pebble $i$.  Initially, all $2k$ pebbles are free.  The goal of Duplicator is
to ensure that every such $\eta_S$ is a partial isomorphism.  If the game reaches an even state $S$
such that $\eta_S$ is not a partial isomorphism, Spoiler wins; otherwise the game continues forever
and Duplicator wins.

\begin{Fact}[\cite{Barwise,bounds}]\label{games fa}
Let $l\leq k$ and consider the version of $\Ga_\tau^k$ where the initial
state is as follows: pebbles $1,\ldots,l$ are placed at elements
$x_1,\ldots,x_l$ of $A$ and at elements $y_1,\ldots,y_l$ of $B$.  If
Duplicator has a winning strategy in that game then, for every $\tau$-formula
$\phi(v_1,\ldots,v_l)\in \Lk$,
\[ A\models\phi(x_1,\ldots,x_l)\; \iff\; B\models\phi(y_1,\ldots,y_l). \]
\end{Fact}

\section{The Deterministic Construction}\label{deterministic se}

We are now ready to state our main theorem:

\begin{Theorem}\label{main th}
There exists a set of finite directed graphs, $\G = \set{G_1,G_2,\ldots}$,
such that $\G$ admits fixed points of unbounded depth and yet on $\G$, $\fo
= (\fo +\lfp)$, i.e.  every formula expressible with a least fixed point
operator is already first-order expressible.
\end{Theorem}

The proof of Theorem \ref{main th} has two main ideas.  The first is
the idea of a standard oracle construction from
Structural Complexity Theory.  The second is Lemma
\ref{clique le}:  a formula in $(\fo +\lfp)$ with only $k$
distinct variables cannot distinguish a $k$-clique from
any larger clique.  We divide the proof up into several
parts, that of the oracle construction (Section 3.1),
that with one free variable (Section 3.2), and finally
the general case (Section 3.3).

\subsection{With Lots of Relation Symbols}

In this subsection we concentrate on
the oracle construction by temporarily introducing infinitely many
new relation symbols of each arity:  $R_i^j$, $i,j\geq
1$.  For convenience in the proofs we will use the
notation $\var(\phi)$ to denote the number of distinct
variables free or bound occurring in $\phi$.  Let
$\free(\phi)$ denote the number of free variables
occurring in $\phi$.

\begin{Lemma}\label{S le}
There exists a set of finite directed graphs, $\D =
\set{D_1,D_2,\ldots}$, which also interpret the new
relations: $R_i^j$, $i,j\geq
1$,
such that $\D$ admits fixed points of unbounded depth; and yet on $\D$, $\fo
= (\fo +\lfp)$, i.e., every formula expressible with a least fixed point
operator is already first-order expressible.
\end{Lemma}

\begin{proof}
Let $\Delta_1,\Delta_2,\ldots$ be a listing of all
formulas in $(\fo + \lfp)$ in this expanded language.
Let $u_i = \free(\Delta_i)$, the number of free variables
occurring in $\Delta_i$.
Let $S_i$ be one of the new relation symbols
of arity $u_i$ such that,
\begin{myequation}{proviso eq}
{\hbox{$S_i$ {\rm does not occur in $\Delta_r$ for
any $r\leq i$.}}}
\end{myequation}

We will let the graph $D_j^0=\angle{V_j,E_j}$ be a directed segment of
length $j-1$:
\begin{eqnarray*}
V_j &=&\set{d_1,d_2,\ldots,d_j}\\
 E_j &=&\bigset{\angle{d_k,d_{k+1}}}{0\leq k<j}
\end{eqnarray*}

We next show how to interpret the new relation symbols in
the $D_j$'s such that: For all $i$, for all $j\geq i$,
and for all $a_1,a_2,\ldots,a_{u_i}\in\abs{D_j}$,
\begin{myequation}{S eq}
{D_j \models (\Delta_i(a_1,a_2,\ldots,a_{u_i})
\,\leftrightarrow\, S_i(a_1,a_2,\ldots,a_{u_i}))}
\end{myequation}

\relax FromEquation \ref{S eq}, it follows that each $\Delta_i$
is equivalent to a first-order formula -- in fact, to an
atomic formula -- for all but finitely many structures.  
Of course, on any fixed finite structure, the formula $\Delta_i$ is
equivalent to a first-order formula.
Lemma \ref{S le} follows immediately.

Now we construct the $D_j$'s so that Equation
\ref{S eq} holds.  $D_j^0$ defined above is just a graph,
which may be thought of as interpreting all of the new relations
as false.  Assuming $D_j^{i-1}$ has been defined, let
$D_j^i$ be the same as $D_j^{i-1}$ except that for all
$a_1,a_2,\ldots,a_{u_i}\in\abs{D_j}$, we interpret $S_i$
so that
\[D_j \models (\Delta_i(a_1,a_2,\ldots,a_{u_i})
\,\leftrightarrow\, S_i(a_1,a_2,\ldots,a_{u_i}))\]
Note that by Equation \ref{proviso eq}, this doesn't affect any of the previous
steps.  

Let $D_j = D_j^j$.  This completes the construction,
guaranteeing that Equation \ref{S eq} holds.  This
completes the proof of Lemma \ref{S le}.
\end{proof}

\subsection{One Free Variable Case: Relations Replaced by Cliques}

Now, we get rid of the new relation symbols, replacing
them by cliques attached to the vertices in the $D_j$'s.
The main result we will need is
that formulas from $\Lk$, i.e. infinitary formulas with at most $k$
variables, cannot distinguish $k$-cliques from $r$-cliques
for any $r>k$.

\begin{Lemma}\label{clique le}
Let $F$ be a finite, directed graph and let $v$ be a vertex in
$F$.  For $i\geq 1$,
let  $F_i$ be the result of replacing $v$ by a clique of
$i$ new vertices: $v_1,\ldots, v_i$.  Each edge
$\angle{v,w}$ or $\angle{z,v}$ from $F$ is replaced with
$i$ new edges: $\angle{v_j,w}$ or $\angle{z,v_j}$, $j=1,2,\ldots,i$.
Let $1\leq k<r$ be natural numbers.  Then 
$F_k$ and $F_r$ agree on all formulas with at most $k$
variables from $\L$.
\end{Lemma}

\begin{proof}
This is proved by using the game $\Ga_\tau^k$ from Fact \ref{games fa}.  We have to show
that the Duplicator has a winning strategy for the $k$-pebble game on $F_k$ and $F_r$.  Her strategy
is to answer any move outside of the cliques with the same vertex in the other graph.  A move on one
of the new cliques is likewise matched by a move on the new clique in the other graph.  Since there
are only $k$ pebbles, there is always an unpebbled vertex in either of the cliques to match with.
Thus the Duplicator has a winning strategy.  It follows that $F_k$ and $F_r$ agree on all
formulas from $\Lk$.
\end{proof}

To make the deterministic construction easier to understand we begin by doing it just for formulas
with only one free variable:

\begin{Lemma}\label{H le}
There exists a set of finite directed graphs, $\H =
\set{H_1,H_2,\ldots}$, such that $\H$ admits fixed points
of unbounded depth, and yet on $\H$, every formula with at
most one free variable that is expressible with a least fixed point
operator is already first-order expressible.
\end{Lemma}

\pf
Let $\Theta_1,\Theta_2,\ldots $ be the set of all formulas
in $(\fo + \lfp)$ that have at most one free variable.  
The construction of the $H_j$'s is similar to that of the
$D_j$'s of Lemma \ref{S le}.  The difference is that
instead of making the relation $S_i(d)$
hold, we will modify the size of a certain clique that is
connected to $d$.  

We next define the sequence of natural numbers: $v_0 <
v_1 < v_2 <\cdots$ that will be the sizes of the initial
cliques.  Let $v_0=0$, and inductively, let $v_{i}
=\max(\var(\Theta_{i}),v_{i-1} + 2^{i+1})$.  In the
construction of $H_j$ we will modify the sizes of cliques
that are initially of size $v_i$.  The modification will
add a number of vertices to these cliques while keeping
them smaller than $v_{i+1}$.

Define the graph $H_j^0$ as follows:  First, $H_j^0$
contains $D^0_j$, the directed segment of length $j-1$.
For each $d\in \abs{D^0_j}$ and for each $i\leq j$,
$H_j^0$ also contains the size $v_i$ clique $C_{d,i}$
which has edges from each of its elements to the vertex
$d$.

Assuming $H_j^{i-1}$ has been defined, let $H^i_j$ be the
same as $H_j^{i-1}$ except that for each $d\in
\abs{D_j^0}$ we add $n(d,i)$ vertices to the $v_i$-vertex clique
$C_{d,i}$.  The number $n(d,i)$ is an $i+1$ bit binary
number such that:
\[(\mbox{``Bit 0 of $n(d,i)$ is one.''}) \quad\Leftrightarrow\quad
( H_j^{i-1} \models \Theta_i(d)) \]
And, for $1\leq s\leq i$, let $a_s$ be a vertex in $C_{d,s}$.
Then,
\[(\mbox{``Bit s of $n(d,i)$ is one.''}) \quad\Leftrightarrow\quad
( H_j^{i-1} \models \Theta_i(a_s)) \]

Finally, let $H_j = H_j^j$.
Define the notation $S\esub_k T$ to mean that $S$ is a
$k$-variable elementary substructure of $T$.  That
is, $S$ is a substructure of $T$ and for all first-order
formulas $\phi$ with $\var(\phi)\leq k$, and for all
$a_1,a_2,\ldots,a_k\in \abs{S}$,
\[ S\models \phi(a_1,a_2,\ldots,a_k)
\quad\Leftrightarrow\quad
  T\models \phi(a_1,a_2,\ldots,a_k)\]
We have constructed the $H_j$'s so that,

\begin{myequation}
{esub eq}
{H_j^{i-1} \; \esub_{v_i} \; H_j}
\end{myequation}

Equation \ref{esub eq} follows from Lemma \ref{clique le}
and the fact the the construction of $H_j^r$ for $r\geq
i$ proceeds by increasing the size of cliques whose size
is at least $v_i$.

Let $a\in\abs{H_j}$.  If $a = d\in\abs{D_j^0}$ then,
\begin{eqnarray*}
(H_j\models \Theta_i(a)) &\Leftrightarrow &
(H_j^{i-1}\models\Theta_i(d))\\
&\Leftrightarrow &
(\mbox{``Bit 0 of $n(d,i)$ is one.''})
\end{eqnarray*}

If $a$ is a member of a clique $C_{d,r}$, let
$s=\min(i,r)$.  Then,
\begin{eqnarray*}
(H_j\models \Theta_i(a)) &\Leftrightarrow &
(H_j^{i-1}\models\Theta_i(a))\\
&\Leftrightarrow &
(\mbox{``Bit s of $n(d,i)$ is one.''})
\end{eqnarray*}

Remember that $v_{i+1}$ is a fixed constant.  Furthermore,
there are at most $2^{i+1}$ possible values for $n(d,i)$.
It follows that there is a first-order formula
$\phi_i(a)$ that finds the appropriate $d$ and $s$,
and determines $n(d,i)$ which is the size of largest
maximal clique connected to $d$ that has fewer than
$v_{i+1}$ vertices.  Next, compute bit $s$ of
$n(d,i)$ by table look up, and let $\phi_i(a)$ be true
iff this bit is one.

Thus, we have that for all $j\geq i$ and for all
$a\in\abs{H_j}$,

\smallskip\noindent\makebox[3.25in]{\hfill${\displaystyle H_j 
\models (\Theta_i(a) \;\leftrightarrow\; \phi_i(a)) }$\hfill $\Box$}

\subsection{General Case: Arbitrary Arity}

The reason that the general case is more complicated than the arity one case is that we must include
gadgets that identify tuples of nodes.  We then must contend with having arguments from these
gadgets and so the arities seem to multiply.  We must therefore be careful so that the arities
remain bounded.

\begin{proof}{\bf of Theorem} \ref{main th}:
Let $\Gamma_1,\Gamma_2,\ldots $ be a listing of all formulas in $(\fo +\lfp)$.
As we have mentioned, arities might multiply.  The base
arity of the formula $\Gamma_i$ is $f_i=\free(\Gamma_i)$.
We will use increased arities $A_0<A_1 <\ldots< A_j$
defined by $A_0=1$, and inductively,
\begin{myequation}
{A eq}
{ A_i \; =\; 1 + (A_{i-1})(2f_i)}
\end{myequation}

Next define the sequence of natural numbers: $w_0 <
w_1 < w_2 <\cdots$ that will be the sizes of the initial
cliques.  Let $w_0=0$, and inductively, let $w_{i} =
\max(\var(\Gamma_{i}),1 + w_{i-1}+ A_{i-1})$.

To define the graph $G_j$, we begin as usual by including
the directed segment $D_j^0$.  For each $i$, we include
enough gadgets: $T_i^r$, $r =1,2,\ldots,n_i$, to encode
all possible sequences of length at most $A_i$ of
elements of $\abs{D_j^0}$.  (Here, $n_i$ is equal to $(j+1)^{A_i}$.)

Each gadget $T_i^r$ consists of $j\cdot A_i$ cliques of
size $w_i$.  For each $d\in\abs{D_j^0}$ there are $A_i$
of these cliques, $C^r_{d,i}$, with edges to $d$.
$T_i^r$ also contains one vertex $t^r_i$
with edges to all the $C^r_{d,i}$'s, $d=1,\ldots,j$.
When we
want $T_i^r$ to encode the sequence
$d_1,d_2,\ldots,d_{A_i}$ we will choose $A_i$ cliques,
$C^r_{d_1,i},C^r_{d_2,i},\ldots,C^r_{d_{A_i},i}$ and
increase their sizes by $1,2,\ldots, A_i$ vertices
respectively.  Note that we have enough copies of each
$C^r_{d,i}$ to tolerate any number of repetitions of the
same $d$.  To skip one of the members of the sequence,
say $d_t$, we increase no clique by exactly $t$ vertices.
In this case we write $d_t = 0$.  Thus, we have
shown how to modify the gadget $T^r_i$ so that it codes
any sequence of length $A_i$ from the alphabet
$\set{0,1,\ldots,j}$.  Note that no formula
$\Gamma_t$ with $t\leq i$ can detect this modification!

Define $G_j^0$ to include $D_j^0$ plus all of the
$T_i^r$'s, $1\leq i\leq j$, $1\leq r\leq n_i$.

Inductively, assume that $G_j^{i-1}$ has been constructed.
Now, for each tuple $a_1,a_2,\ldots,
a_{f_i}\in\abs{G_j^{i-1}}$, {\bf if }
$ G_j^{i-1}\models \Gamma_i(a_1,a_2,\ldots,a_{f_i})$,
{\bf then }
 we will modify one of the
gadgets $T_i^r$ to encode the tuple, $a_1,a_2,\ldots,
a_{f_i}$.

Let's first consider the case that $a_1$ is a vertex from
some $T^{r_1}_{i-1}$.  In this case, $T^{r_1}_{i-1}$ codes a
sequence,
\begin{myequation}
{subseq eq}
{b_{11},b_{12},\ldots,b_{1,A_{i-1}},\;\mbox{
each }b_{1t}\in\set{0,1,\ldots,j}}
\end{myequation}

To reencode this sequence, we first just copy it.  Next,
we have to indicate which vertex in $T^{r_1}_{i-1}$, $a_1$
is.  (It could be the vertex $t^{r_1}_{i-1}$, or a vertex in one of the
unused cliques, $C^{r_1}_{d,i-1}$, or in one of the
cliques $C^{r_1}_{b_{1q},i-1}$ that codes the $q\th$
element of the sequence of Equation \ref{subseq eq}.  In
each case, we use the $A_{i-1}$ extra slots to encode
which of these cases apply\footnote{For those who
want to know, the coding is done as follows: If $a_1$ is
the vertex $t_{i-1}^{r_1}$, then the extra $A_{i-1}$
slots are all 0's.  If $a_1$ is in an unused
$C^{r_1}_{d,i-1}$, then the first two extra slots contain
$d$'s and the rest are 0's.  Finally, if $a_1$ is in
$C^{r_1}_{b_{1q},i-1}$ then put $b_{1q}$ into the $q\th$
extra slot and leave the rest 0.}.
This is the reason for the
factor of 2 in Equation \ref{A eq} and while this is
slightly wasteful, it is simple and we are just trying to
prove that something is finite.

We have just explained how to encode $a_1$ in the first
$2A_{i-1}$ slots of $T_i^r$.  Similarly, code
$a_2,\ldots,a_{f_i}$ into the next $2A_{i-1}(f_i -1)$
slots.  (If one of the $a_s$'s comes from a shorter
sequence, then leave the rest of its positions 0.)
Finally, in the one remaining slot, put a 1.

Let $G_j = G_j^j$.  It follows just as in
Equation \ref{esub eq} that, $G_j^{i-1} \esub_{w_i}G_j$.

Again recall that each $A_i$ and $w_{i+1}$ is a fixed
constant.  Thus, given a tuple, $a_1,\ldots,a_{f_i}$ from
$\abs{G_j}$, a first-order formula,
$\psi_i(a_1,\ldots,a_{f_i})$, can express the existence of
the gadget
$T_i^r$ that codes this tuple.
Thus, for all $j\geq i$,
\[G_j\,\models (\Gamma_i(a_1,\ldots,a_{f_i}) \;\leftrightarrow\,
\psi_i(a_1,\ldots,a_{f_i}))\]
This complete the proof of Theorem \ref{main th}.
\end{proof}

We should note that Theorem \ref{main th} did not use any
properties of $(\fo + \lfp)$ except that the language is
countable and each formula had a constant number of
variables.  We thus have the following extension:

\begin{Corollary}\label{L co}
Let $\sel$ be any countable subset of formulas about
graphs from $L^\omega_{\infty\omega}$.  Then there
exists a set of finite graphs, ${\cal F}$, that admits
unbounded fixed points and such that over
${\cal F}$
every formula from $\sel$ is equivalent to a first-order formula.
\end{Corollary}

\subsection{Two Extensions and an Open Problem}

The deterministic construction relied heavily
on Lemma \ref{clique le}.  This in turn depends on
the fact that $\L$ on unordered structures is
not expressive enough to count.

In \cite{CFI} a lower bound was proved on the language
$(\fo + \cownt + \lfp)$.  This is a language over
two-sorted structures: one sort is the numbers:
$\set{0,1,\ldots, n-1}$ equipped with the usual ordering.
The other sort is the vertices: $\set{v_0,v_1,\ldots,
v_{n-1}}$ with the edge predicate.  The interaction
between the two sorts is via counting quantifiers.  For
example, the formula,
\[\cq {i}{x}\phi(x)\]
means that there exist at least $i$ vertices $x$ such that $\phi(x)$.  Here $i$ ranges over numbers
and $x$ over vertices.  The least fixed point operator may be applied to relations over a
combination of number and vertex variables.  Define the language $(L
+\cownt)_{\infty,\omega}^\omega$ to be the superset of $(\fo +
\cownt + \lfp)$ obtained by adding counting quantifiers to $\L$.

In \cite{CFI} it is shown that the language  $(\fo +
\cownt + \lfp)$ -- and in fact even $(L +\cownt)_{\infty,\omega}^\omega$ -- does not express all polynomial-time properties,
even over structures of color class size four.  Such structures are ``almost ordered'': they consist
of an ordered set of $n/4$ color classes, each of size four.  Only the vertices inside these color
classes are not ordered.  We glean the following fact from \cite{CFI}.

\begin{Fact}[\cite{CFI}]\label{cfi fa}
For each $n>0$ there exist nonisomorphic graphs $T_n$ and $\widetilde{T_n}$ each with $O(n)$
vertices, such that $T_n$ and $\widetilde{T_n}$ are indistinguishable by all formulas with at most
$n$ variables from $(\fo +\lfp + \cownt)$, or even from $(L +\cownt)_{\infty,\omega}^\omega$.
\end{Fact}

Useful in the proof of Fact \ref{cfi fa} as well as in the next theorem is the 
following modification of the game $\Gamma_\tau^k$ of Fact \ref{games fa}.
Given a
pair of $\tau$-structures $G$ and $H$ define the $\C_\tau^k$ game on $G$ and $H$ as
follows: Just as in the $\Gamma_\tau^k$ game, we have two players and
$k$ pairs of pebbles.  The difference is that each move now
has two parts.
\begin{enumerate}
\item Spoiler picks up the pair of pebbles numbered $i$ for some $i$.
He then chooses a set $A$ of vertices from one of the
graphs.  Now Duplicator answers with a set $B$ of vertices
from the other graph.  $B$ must have the same cardinality as
$A$.
\item Spoiler places one of the pebbles numbered $i$ 
on some vertex $b\in B$.  Duplicator
answers by placing the other pebble numbered $i$ on some $a\in A$.
\end{enumerate}

The definition for winning is as before.  What is going on
in the two part move is Spoiler asserts that there
exist $\abs{A}$ vertices in $G$ with a certain property.
Duplicator answers with the same number of such vertices in
$H$.  Spoiler challenges one of the vertices in $B$ and
Duplicator replies with an equivalent vertex from $A$.  This
game captures expressibility in $(L +\cownt)_{\infty,\omega}^\omega$:

\begin{Fact}[\cite{canon}]
The Duplicator has a winning strategy for the $\C_\tau^k$ game on $G,H$ if and only if $G$ and $H$
agree on all formulas with at most $k$ variables from $(L +\cownt)_{\infty,\omega}^\omega$.
\end{Fact}

Using the above facts, we now prove a counterexample to a
weaker version of McColm's Conjecture:

\begin{Theorem}\label{count th}
There exists a set of finite directed graphs, $\J = \set{J_1,J_2,\ldots}$,
such that $\J$ admits fixed points of unbounded depth and yet on $\J$, $\fo
= (\fo+\cownt +\lfp)$, i.e., every formula expressible with a least fixed point
operator and counting is already first-order expressible.  In fact, this statement remains true
when $(\fo+\cownt +\lfp)$ is replaced by an arbitrary countable subset of $(L +\cownt)_{\infty,\omega}^\omega$.
\end{Theorem}

\begin{proof}
The idea of this construction is that everywhere we
started with
a clique of size $n$ in the previous proof, we will
start with
a chain of copies of the graph $T_n$ from Fact \ref{cfi
fa}.  Then where previously we increased the size of the
clique to code some number $b$ of bits, we will instead flip
some copies of $T_n$ to
$\widetilde{T_n}$, in a particular length $b$ chain of $T_n$'s.  

The main differences are that unlike
the cliques, there is not an automorphism mapping every
point in $T_n$ to every other point in $T_n$.  Furthermore, $T_n$ is distinguishable from $T_{n+1}$
using a small number of variables.  

Let $f(j)$ be the number of formulas that are handled by the structure $G_j$, and let $v(j)$ be
$v_{f(j)}$, the number of variables to be handled as in the proof of Theorem \ref{main th}.  Observe
that $f(j)$ and thus $v(j)$ may be chosen to grow very slowly.  In particular, we will make sure
that $f(j)$, and in fact the number of vertices in each $T_{v(j)}$ is less than $j$.  Recall also
that the graphs $T_n$ from Fact \ref{cfi fa} are ordered up to sets of size four.  We introduce two
new binary relations: Red edges from each vertex in each $T_{v(i)}$ to the vertex $i\in D_j^0$, and
Blue edges from each of the four vertices numbered $k$ in any of the $T_{v(i)}$'s to the vertex
$k\in D_j^0$.  Thus, any vertex chosen from $G_j$ will have a ``name'' that consists of a pair of
vertices from $D_j^0$, together with a  bounded number of bits.  

The construction and proof now follow as in the proof of Theorem \ref{main th}.
\end{proof}

We also show,

\begin{Corollary}\label{ppspace co}
If $\p\ne\pspace$, then there exists a set $\C$ of finite structures such that $\fo = (\fo +\lfp)$
on $\C$; but, $\fo \ne (\fo + \iterp)$ on $\C$.
\end{Corollary}

\begin{proof}
Let $\G$ be the set of all finite, ordered graphs.  If $\p\ne\pspace$, then there is a property
$S\subset \G$ such that $S\in \pspace - \p$.  Now, do the construction of Theorem \ref{main th},
starting with $\G$.  This construction assures that $\fo = (\fo +\lfp)$
on the resulting set $\C$.  However, any first-order formula $\phi$ has a fixed number, $k$, of
variables.  Thus, to $\phi$, the noticeable changes during the construction involve at most $k$
PTIME properties.  Therefore, $S$ is still not recognizable in $\fo$ over $\C$.
\end{proof}

One special case of McColm's conjecture remains open.  This
is a fascinating question in complexity theory and logic
related to uniformity of circuits and logical
descriptions, cf. \cite{BIS}.  Consider the structures
$\B = \set{B_1,B_2,\ldots}$ where $B_i =
\angle{\set{0,1,\ldots,i-1},\leq,\bit}$.
Here $\leq$ is the usual ordering on the natural numbers
and $\bit(x,y)$ holds iff the $x\th$ bit in the binary
representation of the number $y$ is a one.

\begin{Question}\label{uniform qu}
Is $\fo = (\fo + \lfp)$ over $\B$?
\end{Question}

The answer to Question\ \ref{uniform qu} is ``Yes,'' iff every polynomial-time computable numeric
predicate is already computable in $(\fo + \bit)$.  Equivalently, the answer to Question
\ref{uniform qu} is ``Yes,'' iff deterministic logtime uniform $\ac^0$ is equal to polynomial-time
uniform $\ac^0$, cf. \cite{BIS}.  A resolution of this question would thus answer an important
question in complexity theory.

\section{The Randomized Construction}\label{random se}

We now sketch a quite different construction that also disproves McColm's conjecture.  
Throughout this construction, $P$ is a binary predicate.  We will prove:

\begin{Theorem}\label{TheoremOne}
Suppose that $K_1$ is a class of structures of some vocabulary $\tau_1$, and
$\cL$ is an arbitrary countable subset of $\L$.  Let $\tau_2$ be the
extension of $\tau_1$ with an additional binary predicate $P$.  There exist
a class $K_2$ of $\tau_2$-structures such that:
\begin{enumerate}\item
$K_1$ is precisely the class of $\tau_1$-reducts of substructures
$M_2\mid\{x\mid P(x,x)\}$ where $M_2$ ranges over $K_2$.
\item
Every $\cL$-formula is equivalent to a first-order formula in $K_2$.
\end{enumerate}
\end{Theorem}

The idea of the proof is relatively simple.  Let $\rho_1,\rho_2,\ldots$ be a
list of all $\cL$-definable global relations on $K_1$.  We attach a graph
$G$ to every $M\in K_1$ and define a projection function from elements of
the new sort to elements of the old sort.  Relations $\rho_i^M$ on the old
sort are coded by cliques of $G$; a tuple $\a$ belongs to $\rho_i^M$ if and
only if there is clique of cardinality $i$ projected in a certain way onto
$\a$.  The necessity to have appropriate cliques is the only constraint on
$G$; otherwise the graph is random.  We check that every $\cL$-definable
global relation reduces by first-order means to $\cL$-definable global
relations on the old sort and thus is first-order expressible.  In fact, we
beef $\cL$ up before executing the idea.

Let $H$ be a hypergraph of cardinality $\geq 2$.

\begin{Definition}\label{EnvelopeDefinition}
An {\em envelope} for $H$ is a $\{P\}$-structure $E$ satisfying the 
following conditions:
\begin{itemize}\item
$|H|\subseteq|E|$, and $P$ is the identity relation on $|H|$.
\item
$P$ is irreflexive and symmetric on $|E|-|H|$.
\item
For every $x\in|E|-|H|$, there is a unique $a\in H$ with $E\models P(x,a)$.
\item
For every $a\in|H|$ and every $x\in|E|-|H|$, $E\models\neg P(a,x)$. 
\end{itemize}
\end{Definition}

Let $E$ range over envelopes for $H$ such that $|E|-|H|\neq\emptyset$.

\begin{Definition}
Elements of $H$ are {\em nodes\/} of $E$ and elements of $|E|-|H|$ are {\em
vertices\/} of $E$.  $G_E$ is the graph formed by $P$ on the vertices.  If
$E\models P(x,a)$ and $a\in H$ then $a$ is called the {\em projection\/} of
$x$ and denoted $F(x)$ (or $Fx$).  If $X$ is a set of elements of $E$ then
$F(X)$ is the multiset
\( \openmultiset Fx \mid x\in X \closemultiset \). 
If $\x$ is a sequence $(x_1,\ldots,x_l)$ of elements of $E$ then
\( F(\x) = ( F(x_1),\ldots,F(x_l) ) \).
\end{Definition}

Let $k$ be a positive integer $\geq 3$.

\begin{Definition}
A clique $X$ of $G_E$ is a {\em $k$-clique\/} if $F(X)\in\HE(H)$ and
$||X||<k$.  A vertex that does not belong to any $k$-clique is {\em
$k$-plebeian}.  The {\em $k$-closure\/} $C_k(X)$ of a subset $X$ of $E$ is
the union of $X$ and all $k$-cliques intersected by $X$.
\end{Definition}

\begin{Definition}
$E$ is {\em $k$-good\/} for $H$ if it satisfies the following conditions.
\begin{description}\item[$G_0(k)$ ]
All $k$-cliques are pairwise disjoint.

\item[$G_1(k)$ ]
For every $X\subseteq|E|$ of cardinality $<k$, there is a $k$-plebeian
vertex $z\in|E|-X$ with a predefined projection $Fz$ which is $P$-related to
$C_k(X)$ in any predefined way that does not destroy any $k$-clique
$C\subseteq C_k(X)$.  In other words, if $a$ is a node, $Y\subseteq C_k(X)$
and $Y$ does not include any $k$-clique, then there is a $k$-plebeian vertex
$z\in F^{-1}(a)-X$ adjacent to every vertex in $Y$ and to no vertex in
$C_k(X)-Y$.

\item[$G_2(k)$ ]
For every $X\subseteq|E|$ of cardinality $<k$, there is a $k$-clique
$\{y_1,\ldots,y_l\}\subseteq|E|-X$ with any predefined projections $Fy_m$
and any predefined pattern
\( R = \{ (x,m) \mid E\models P(x,y_m)\} \)
that does not destroy any $k$-clique $C\subseteq C_k(X)$.  In other words,
if $\a=(a_1,\ldots,a_l)$ is a tuple of nodes, $l<k$, $\MS(\a)$ is a
hyperedge, $R\subseteq C_k(X)\times\{1,\ldots,l\}$, no vertex is
$R$-adjacent to all the numbers, and no number is $R$-adjacent to all
vertices of any $k$-clique $C\subseteq C_k(X)$, then there is a tuple
$\y=(y_1,\ldots,y_l)$ of distinct vertices such that $F(\y)=\a$,
$\{y_1,\ldots,y_l\}$ is a clique disjoint from $X$, and $E\models P(x,y_m)
\iff (x,m)\in R$ for all $x\in C_k(X)$ and all $m$.
\end{description}
\end{Definition}

\begin{Lemma}\label{GoodEnvelopeLemma}
\begin{enumerate}\item
If $E$ is $k$-good, $X\subseteq E$ and $||E||<k$ then
$||C_k(X)||\leq(k_1)^2$. 
\item
If $E$ is $k$-good then every hyperedge of cardinality $<k$ is the
projection of some $k$-clique.
\item
In every $k$-good envelope, every clique $C$ of cardinality $<k$ is a
$k$-clique.  Moreover, if a clique $C\subseteq C_k(X)$ for some $X$ of
cardinality $<k$ then $C$ is a $k$-clique.
\item
Let $H'$ be the hypergraph obtained from $H$ by discarding all hyperedges of
cardinality $\geq k$.  Then $E$ is $k$-good for $H$ if and only if it is
$k$-good for $H'$.
\item
If $E$ is $k'$-good for $H$ where $k'>k$ then $E$ is $k$-good for $H$.
\end{enumerate}
\end{Lemma}

\begin{proof}
Omitted due to lack of space.
\end{proof}

\begin{Theorem}\label{GoodEnvelopeTheorem}
There exists a $k$-good envelope for $H$.
\end{Theorem}

\begin{proof}
Omitted due to lack of space.
\end{proof}

\subsection{The Game}

Let $M$ be a structure of some vocabulary $\tau_0$ such that every element
of $M$ interprets some individual constant.  It is supposed that $\tau_0$
does not contain the fixed binary predicate $P$.  Let $H$ be a hypergraph on
$|M|$, so that $|H|=|M|$.  An envelope $E$ for $H$ can be seen as a
structure of vocabulary $\tau=\tau_0\cup\{P\}$ where the $\tau_0$-reduct of
the substructure $E\mid|H|$ equals $M$ and no $\tau_0$ relation involves
elements of $|E|-|H|$.

Fix a positive integer $k$ and let $E$ and $E'$ range over $k$-good
envelopes for $H$.  We will prove that Duplicator has a winning strategy in
$\Ga_\tau^k(E,E')$.

\begin{Definition} 
A partial isomorphism $\eta$ from $E$ to $E'$ is {\em $k$-correct\/} if it
satisfies the following conditions where $x$ ranges over $\Dom(\eta)$. 
\begin{itemize}\item
If $x$ is a node then $\eta(x)=x$.
\item
If $x$ is a vertex then $\eta(x)$ is a vertex and $F(\eta(x))=Fx$.
\item 
$x$ is $k$-plebeian if and only if $\eta(x)$ is $k$-plebeian.
\item
If $x$ belongs to some $k$-clique $X$ then $\eta(x)$ belongs to some
$k$-clique $X'$ such that $F(X')=F(X)$.
\end{itemize}
\end{Definition}

\begin{Definition}
A $k$-correct partial isomorphism $\eta$ from $E$ to $E'$ is {\em $k$-nice}
if there exists an extension of $\eta$ to a $k$-correct partial isomorphism
$\eta^*$ with domain $C_k(\Dom(\eta))$.
\end{Definition}

\begin{Lemma}\label{EtaStarLemma}
Suppose that $\eta$ is a $k$-nice partial isomorphism from $E$ to $E'$.
Then $\eta^*$ and $\eta^{-1}$ are $k$-nice, $(\eta^*)^{-1}=(\eta^{-1})^*$,
and $\Range(\eta^*)=C_k(\Range(\eta))$.  $\eta^*$ maps every $k$-clique onto
$k$-clique of the same size, different $k$-cliques are mapped to different
$k$-cliques.
\end{Lemma}

\begin{proof}
Obvious.
\end{proof}

\begin{Definition}
An even-numbered state of $\Ga_\tau^k(E,E')$ is {\em good\/} if the
pebble-defined map is a $k$-nice partial isomorphism.  A strategy of
Duplicator in $\Ga_\tau^k(E,E')$ is {\em good\/} if every move of Duplicator
creates a good state.
\end{Definition}

\begin{Theorem}\label{good strategy th}
Every good strategy of Duplicator wins $\Ga_\tau^k(E,E')$, and Duplicator
has a good strategy.
\end{Theorem}

\begin{proof}
Omitted due to lack of space.
\end{proof}

\begin{Definition}
A {\em $0$-table\/} is a conjunction $\al(v_1,\ldots,v_l)$ of atomic and
negated atomic formulas in vocabulary $\{P\}$ which describes the
isomorphism type of a $\{P\}$-structure of cardinality $\leq l$ which can be
embedded into some envelope for some hypergraph.
\end{Definition}

\begin{Definition}
Let $j<k$ be a positive integer.  A {\em $(j,k)$-table\/} is a first-order
$\{P\}$-formula $\be(v_1,\ldots,v_l)$ which says that there are distinct
elements $u_1,\ldots,u_j$ such that $\{u_1,\ldots,u_j\}$ is a clique
intersecting $\{v_1,\ldots,v_l\}$ and a particular $0$-table
$\be_0(u_1,\ldots,u_s,v_1,\ldots,v_l)$ is satisfied.
\end{Definition}

\begin{Definition}
A {\em $k$-table\/} $\ga(v_1,\ldots,v_l)$ is a conjunction such that:
\begin{itemize}\item
Some $0$-table $\al(v_1,\ldots,v_l)$ is a conjunct of $\ga(v_1,\ldots,v_l)$.
\item
If $j<k$ and $\be(v_1,\ldots,v_l)$ is a $(j,k)$-table consistent with
$\al(v_1,\ldots,v_l)$ then either $\be(v_1,\ldots,v_l)$ or
$\neg\be(v_1,\ldots,v_l)$ is a conjunct of $\ga(v_1,\ldots,v_l)$.
\item
There are no other conjuncts.
\end{itemize}
\end{Definition}

Fix a $k$-variable infinitary $\tau$-formula
$\phi(u_1,\ldots,u_l,v_1,\ldots,v_m)$ and let $\Phi(\u,\v)$ be the
conjunction of $\phi(\u,\v)$ and some $k$-table $\ga(\v)$.  Let $\a$ be an
$l$-tuple of nodes of $H$ and $b$ be an $m$-tuple of nodes $H$.  We
introduce a relation $\Phi^-(\u,\v)$ on $H$.

\begin{Definition}\label{PhiMinusDefinition}
\[\Phi^-(\a,\b) \iff E\models(\exists\v)[(\Phi(\a,\v))\wedge F(\v)=\b]. \]
\end{Definition}

\begin{Lemma}\label{PhiMinusLemma}
$\Phi^-$ does not depend on the choice of $E$: any other $k$-good envelope
for $H$ yields the same relation.
\end{Lemma}

\begin{proof}
It suffices to check that $E'$ yields the same relation.  Since Duplicator
has a winning strategy in $\Ga_\tau^k(E,E')$, no infinitary $k$-variable
$\tau$-sentence distinguishes between $E$ and $E'$.  In particular, no
sentence
\[(\exists v_1,\ldots,v_m) [\, P(v_1,d_1)\, \wedge\,
\ldots\,\wedge P(v_m,d_m) \] 
\[ \qquad\qquad\wedge\quad  \Phi(c_1,\ldots,c_l,v_1,\ldots,v_m)], \]
where $c_1,\ldots,c_l,d_1,\ldots,d_m$ are individual constants,
distinguishes between $E$ and $E'$.
\end{proof}

\begin{Theorem}\label{QuantifierEliminationTheorem}
Let $\x$ be an $m$-tuple of vertices in $E$.  The following claims are
equivalent:
\begin{enumerate}\item
$E\models\Phi(\a,\x).$
\item
$H\models\Phi^-(\a,F(\x))$ and $E\models\ga(\x)$.
\end{enumerate}
\end{Theorem}

\begin{proof}
Omitted due to lack of space.
\end{proof}

In the case $m=0$, $\Phi=\Phi^-=\phi$ and we have the following corollary.

\begin{Corollary}\label{QuantifierEliminationCorollary}
\[ E\models\phi(\a) \iff H\models\phi(\a). \]
\end{Corollary}

\subsection{Proof of Theorem~\ref{TheoremOne}}

We start with a couple of auxiliary definitions.  Call an $r$-ary relation
$R$ {\em irreflexive\/} if every tuple in $R$ consists of $r$ distinct
elements.  Call a global relation $\rho$ {\em irreflexive\/} if every local
relation $\rho^M$ is so.

\begin{Lemma}\label{IrreflexiveLemma}
Every global relation $\rho(v_1,\ldots,v_r)$ is a positive boolean
combination of irreflexive global relations definable from $\rho$ in a
quantifier-free way.
\end{Lemma}

\begin{proof} 
Omitted due to lack of space.
\end{proof}

Call a multiset $A$ is {\em oriented\/} if the relation $\MP(a)<\MP(b)$ is a
linear order on $\Set(A)$; let $\OSet(A)$ be the corresponding linearly
ordered set. \ss

Now we are ready to prove Theorem~\ref{TheoremOne}.  Suppose that $K_1$ is a
class of structures of some vocabulary $\tau_1$, and $\tau_2$ is the
extension of $\tau_1$ with binary predicate $P$.  Let $\cL$ be an arbitrary
countable set of $\L$-formulas.

A global relation $\rho$ on a class $K$ is {\em decidable\/} if there exists
an algorithm that, given (the encodings of) a structure $M\in K$ and a tuple
$\a$ of elements of $M$ of appropriate length, decides whether
$M\models\rho(\a)$ or not.  We are interested in a relativized version of
this definition where $K$ is the collection of all structures (that is, all
finite structures) in the vocabulary of $\rho$.  Let 
\[ \Om = \{ (\phi,M,\a,1)\mid \phi\in\cL \wedge M\models\phi(\a)\}
\; \cup\]
\[\qquad\qquad \{ (\phi,M,\a,0)\mid \phi\in\cL \wedge M\not\models\phi(\a)\} \]

\begin{Definition}
A global relation $\rho$ of vocabulary $\tau$ is {\em $\cL$-decidable\/} if
there is an algorithm with oracle $\Om$ that, given a $\tau$-structure $M$
and a tuple $\a$ of elements of $M$ of appropriate length, decides whether
$M\models\rho(\a)$ or not.
\end{Definition}

Every global relation defined by a formula in $\cL$ is $\cL$-decidable, and
there there are only countably many $\cL$-decidable relations.  List
all $\cL$-decidable irreflexive global relations on $K_1$ of positive
arities:
\( \rho_2,\rho_3\,\rho_4,\ldots \), 
and let $r_i$ be the arity of $\rho_i$.  We suppose that
\( r_i(r_i+1)/2 \leq i \). 
Let $M$ range over $K_1$ and $i$ range over positive integers $\geq2$.

For each $M$ and each $i$, let $\si_i^M$ be the collection of oriented
multisets $A$ such that $\OSet(A)\in\rho_i^M$ and $||A||=i$.  Since
\( 1+2+\ldots+r_i=r_i(r_i+1)/2\leq i \),
$\si_i^M$ is empty.  Let $H(M)$ be the hypergraph
\[ \left(|M|\:,\: \bigcup \{ \si_i^M \mid 1\leq i\leq||M|| \} \right). \]

Set $\tau_2=\tau_1\cup\{P\}$ and let $\cE(M)$ be the collection of
$||M||$-good envelopes for $H(M)$ of minimal possible cardinality.  (The
minimal cardinality is not important; we will use only the following two
consequences: (i)~$\cE(M)$ is finite, and (ii)~there is an algorithm that,
given $M$ constructs some $E\in\cE(M)$.)  View envelopes $E\in\cE(M)$ as
$\tau_2$-structures where the $\tau_1$-reduct of the substructure $E\mid|M|$
equals $M$ and no $\tau_1$-relation involves elements of $|E|-|M|$.  For
every $K\subseteq K_1$, let
\( \cE(K)=\bigcup_{M\in K} \cE(M) \). 
Finally, let $K_2 = \cE(K_1)$.  By the definition of envelopes
(Definition~\ref{EnvelopeDefinition}), $K_2$ satisfies requirement 1 of
Theorem~\ref{TheoremOne}.  In order to prove requirement 2, it suffices to
prove that every infinitary formula with $\cL$-decidable global relation is
first-order definable in $K_2$.

For any global relation $\rho(\v)$ on $K_1$, let $\rho^+(\v)$ be the global
relation on $K_2$ such that
\[ E\models\rho^+(\x) \iff M\models\rho(F(\x)) \]
if $M\in K$, $E\in\cE(M)$ and $\x$ is a tuple of elements of $E$ of the
appropriate length. 

\begin{Lemma}\label{DecidableRhoLemma}
If $\rho$ is $\cL$-decidable then $\rho^+$ is first-order definable in
$K_2$.
\end{Lemma}

\begin{proof}
Omitted due to lack of space.
\end{proof}

Now let $\phi$ be an arbitrary infinitary $\tau_2$-formula
whose global relation is $\cL$-decidable.  We prove that $\phi$ is
equivalent to a first-order formula in $K_2$.  Without loss of generality,
\( \phi = \phi(u_1,\ldots,u_l,v_1,\ldots,v_m) \)
and $\phi$ implies 
\[P(u_1,u_1),\ldots, P(u_l,u_l),\neg
P(v_1,\v_1),\ldots,\neg P(v_m,v_m)\]
In other words, variables $u_i$ are
node variables, and variables $v_j$ are vertex variables.

Let $k$ be the total number of variables in $\phi$, $K_1' = \{M \mid
||M||\geq k\}$ and $K_2'=\cE(K_1')$, so that every $E\in K_2'$ is $k$-good.
Since $K_2-K_2'$ is finite, it suffices to prove that $\phi(\u,\v)$ is
equivalent to a first-order formula in $K_2'$.  Let $\ga(\v)$ be an
arbitrary $k$-table.  Since there are only finite many $k$-tables, it
suffices to prove that the formula
\( \Phi(\v) = \phi(\v)\wedge\ga(\v) \)
is equivalent to a first-order formula over $K_2'$.

Define a global relation $\Phi^-$ on $K_1$ as follows:
\[M\models\Phi^-(\a,\b) 
\iff (\exists\x)[(E\models\Phi(\x)) \wedge F(\x)=\a] \]
where $E\in\cE(M)$.  The choice of $E$ does not matter.  Indeed, extend
$\tau_1$ with individual constants for each element of $M$; call the
resulting vocabulary $\tau_0$.  Now apply Lemma~\ref{PhiMinusLemma} with
$H=H(M)$.

\begin{Lemma}
$\Phi^-$ is $\cL$-decidable.
\end{Lemma}

\begin{proof}
Clear.
\end{proof}

It is not quite true that $(\Phi^-)^+$ is the global relation of the
formula $\Phi$ on $K_2'$ but this is close to truth.  By virtue of
Theorem~\ref{QuantifierEliminationTheorem},
\[ \Phi(\u,\v) \iff [(\Phi^-)^+(\u,\v) \wedge \ga(\v) \]
on $K_2'$.  Indeed, consider any $M\in K_1'$.  Extend $\tau_1$ with
individual constants for each element of $M$; call the resulting vocabulary
$\tau_0$.  Now apply Theorem~\ref{QuantifierEliminationTheorem} with
$H=H(M)$.  By Lemma~\ref{DecidableRhoLemma}, $(\Phi^-)^+$ is first-order
definable in $K_2$.  It follows that $\Phi$ is equivalent to a first-order
formula on $K_2'$.

\subsection*{Acknowledgment}
We are grateful to Anuj Dawar, Phokion Kolaitis, Steven Lindell, and Scott Weinstein for
stimulating discussions.

\end{document}